\documentclass[12pt,leqno]{article}
\usepackage[mathscr]{eucal}
\usepackage{amsmath,amssymb,latexsym,theorem,bbm,amsfonts}
\usepackage{verbatim}
\usepackage{color}
\usepackage{appendix}
\usepackage{array,epsfig,theorem,bbm,graphics,booktabs,url}

\newcommand{\CC}{\mathbb{C}}

\newcommand{\NN}{\mathbb{N}}

\newcommand{\RR}{\mathbb{R}}
\newcommand{\TT}{\mathbb{T}}

\newcommand{\tf}{\widetilde{f}}

\newcommand{\tg}{\widetilde{g}}

\newcommand{\tbeta}{\widetilde{\beta}}

\newcommand{\tgamma}{\widetilde{\gamma}}

\newcommand{\bone}{{\boldsymbol{1}}}

\newcommand{\cF}{{\mathcal F}}

\newcommand{\tX}{\widetilde{X}}

\newcommand{\dd}{\mathrm{d}}
\newcommand{\ee}{\mathrm{e}}
\newcommand{\ii}{\mathrm{i}}

\newcommand{\EE}{\operatorname{\mathbb{E}}}

\newcommand{\PP}{\operatorname{\mathbb{P}}}

\newcommand{\var}{\operatorname{Var}}
\newcommand{\cov}{\operatorname{Cov}}

\newcommand{\tS}{\widetilde{S}}

\renewcommand{\leq}{\leqslant}
\renewcommand{\geq}{\geqslant}

\newcommand{\proofend}{\hfill\mbox{$\Box$}}

\numberwithin{equation}{section}

\theoremstyle{change} \theorembodyfont{\em}
\newtheorem{Lem}{Lemma.}[section]
\newtheorem{Thm}[Lem]{Theorem.}
\newtheorem{Pro}[Lem]{Proposition.}

\newtheorem{Def}[Lem]{Definition.}

\theorembodyfont{\rm}
\newtheorem{Rem}[Lem]{Remark.}
\newtheorem{Ex}[Lem]{Example.}

\begin{document}

\begin{center}
 {\bfseries\Large
    Dilatively stable stochastic processes\\[1ex] and aggregate similarity}

\vspace*{3mm}

 {\sc\large
  M\'aty\'as $\text{Barczy}^{*,\diamond}$,
  Peter $\text{Kern}^{**}$,
  Gyula $\text{Pap}^{***}$}

\end{center}

\vskip0.2cm

\noindent
 * Faculty of Informatics, University of Debrecen,
   Pf.~12, H--4010 Debrecen, Hungary.

\noindent
 ** Mathematical Institute, Heinrich-Heine-University D\"usseldorf,
    Universit\"atsstr. 1, D--40225 D\"usseldorf, Germany.

\noindent
 *** Bolyai Institute, University of Szeged,
     Aradi v\'ertan\'uk tere 1, H--6720 Szeged, Hungary.

\noindent e--mails: barczy.matyas@inf.unideb.hu (M. Barczy),
                    Peter.Kern@uni-duesseldorf.de (P. Kern),
                    papgy@math.u-szeged.hu (G. Pap).

\noindent $\diamond$ Corresponding author.

\vskip0.2cm

%\centerline{\sl June 12, 2014.}

\renewcommand{\thefootnote}{}
\vspace*{0.2cm}
\footnote{The research has been supported by the DAAD-M\"OB Research Grant
 No.~55757 partially financed by the German Federal Ministry of Education and
 Research (BMBF).}

\vspace*{-10mm}

\begin{abstract}
 Dilatively stable processes generalize the class of infinitely divisible
  self-similar processes.
 We reformulate and extend the definition of dilative stability introduced by
  Igl\'oi (2008) using characteristic functions.
 We also generalize the concept of aggregate similarity introduced by
  Kaj (2005).
 It turns out that these two notions are essentially the same for infinitely divisible processes.
 Examples of dilatively stable generalized fractional L\'evy processes are given and
  we point out that certain limit processes in aggregation models are dilatively stable.
\end{abstract}

\section{Dilative stability and aggregate similarity}
\label{section_intro}

Self-similarity is a scaling property of stochastic processes.
It was Lamperti's paper \cite{Lam} which called the attention to the
 significance of this property (named there semi-stability).
As a generalization, Igl\'oi \cite[Definition 2.1.3 and Theorem 2.2.1]{Igl} introduced
 a more general scaling property of certain infinitely divisible processes called dilative stability.
Igl\'oi \cite[Examples 2.1.5\,-7]{Igl} already provided some important (non self-similar) dilatively stable
 processes, such as non-Gaussian moving-average fractional L\'evy motions.
Roughly speaking, non-Gaussian fractional L\'evy motions are not self-similar
 but they belong to a wider class of processes, to the class of dilatively
 stable processes, which underlines the importance of dilative stability.
To put dilative stability into a more general context, we note that there are many other
 known examples of stochastic processes that do not fit into the scope of self-similarity, but are
 very natural from the point of view of scaling limit theorems and various types of invariance properties
 in probability theory, besides Igl\'oi \cite{Igl}, see, e.g., Bierm\'e et al. \cite{BieEstKaj},
 Kaj \cite{Kaj} and Pilipaysjaut\.{e} and Surgailis \cite{PilSur}.
However, there is no common framework for these kinds of processes.
The examples given by Igl\'oi \cite{Igl} and the new examples presented in Sections \ref{GFLP} and \ref{lithuania}
 show that the class of dilatively stable processes is not only much wider than that of self-similar
 processes, but it contains several important non self-similar processes.
So the systematic study of dilative stability started by Igl\'oi \cite{Igl} deserves to be continued.

Further, Kaj \cite[Section 3.6]{Kaj} introduced a similar concept called
 aggregate similarity appearing for certain heavy-tailed limit processes in teletraffic models under an
 intermediate growth condition.

This paper has three-fold aims.
First, we reformulate and extend the definition of dilative stability introduced by
 Igl\'oi \cite{Igl} using characteristic functions.
 Next, we generalize the concept of aggregate similarity due to Kaj, pointing out that
 these generalized aggregate similarity and dilative stability are essentially the same for infinitely
 divisible processes.
Finally, we present new examples of dilatively stable processes that are not self-similar.

In what follows, \ $\NN$, \ $\RR_+$, \ $\RR_{++}$ \ and \ $\RR$ \ will denote
 the set of positive integers, non-negative real numbers, positive real
 numbers and real numbers, respectively.
Further, let \ $\TT$ \ be either \ $\RR$, \ $\RR_+$ \ or \ $\RR_{++}$.
All the stochastic processes are defined on a probability space
 \ $(\Omega, \cF, \PP)$.

First, we give both definitions: self-similarity and dilative stability.

\begin{Def}\label{DefSS}
Let \ $\alpha > 0$.
\ A process \ $(X_t)_{t\in\TT}$ \ is called $\alpha$-self-similar if in terms of
 characteristic functions it fulfills the space-time scaling relation
 \begin{align}\label{ssrel}
  \EE\left( \ee^{ \ii \sum_{j=1}^k \theta_j X_{Tt_j}} \right)
  = \EE\left( \ee^{ \ii \,T^\alpha\sum_{j=1}^k \theta_j X_{t_j}} \right)
 \end{align}
 for all \ $T > 0$, \ $k \in \NN$, \ $\theta_1, \ldots, \theta_k \in \RR$,
 \ $t_1, \ldots, t_k \in \TT$.
\end{Def}

If \ $0 \in \TT$, \ then \ $\PP(X_0 = 0) = 1$, \ see, e.g., Samorodnitsky and
 Taqqu \cite[page 312]{SamTaq}.
Note that property \eqref{ssrel} is equivalent to
 \begin{align}\label{SS}
  \forall\, T>0 : \, X(T\,\cdot\,) \overset{\text{\text{fd}}}
  \sim T^\alpha X(\,\cdot\,),
 \end{align}
 where $\overset{\text{\text{fd}}}{\sim}$ denotes that the finite-dimensional
 distributions are the same.
If a stochastically continuous process \ $(X_t)_{t\in\RR_+}$ \ is
 $f$-self-similar for some function \ $f : (0, \infty) \to (0, \infty)$ \ in
 the sense that it satisfies \eqref{ssrel} replacing \ $T^\alpha$ \ by \ $f(T)$
 on the right-hand side of \eqref{ssrel},
 then there exists an \ $\alpha \geq 0$ \ such that \ $f(T) = T^\alpha$,
 \ $T > 0$, \ and in case of \ $\alpha = 0$, \ one has \ $\PP(X_t = X_0) = 1$,
 \ $t \in \RR_+$, \ see Lamperti \cite[Theorem 1]{Lam}.
In the definition of self-similarity the case \ $\alpha = 0$ \ is usually
 excluded, such as in Definition \ref{DefSS}.

Dilative stability is an analogous property of certain infinitely divisible
 processes (all finite-dimensional distributions are infinitely divisible)
 involving a scaling also in the convolution exponent.

For an infinitely divisible process \ $(X_t)_{t\in\TT}$ \ and for \ $k \in \NN$,
 \ $t_1, \ldots, t_k \in \TT$, \ let \ $\psi_{t_1,\ldots,t_k} : \RR^k \to \CC$
 \ denote the characteristic exponent of \ $(X_{t_1}, \ldots, X_{t_k})$,
 \ i.e., the unique continuous function such that \ $\psi_{t_1,\ldots,t_k}(0,\ldots,0)=0$ \ and
 \[
   \EE\left( \ee^{ \ii \sum_{j=1}^k \theta_j X_{t_j}} \right)
   = \ee^{- \psi_{t_1,\ldots,t_k}(\theta_1, \ldots, \theta_k)}
   \quad\text{for all \ } \theta_1, \ldots, \theta_k \in \RR ,
 \]
 see, e.g., (the multi-dimensional version of) Lemma 3.2.11 in Stroock \cite{Str}.

\begin{Def}\label{DefDS}
Let \ $\alpha, \delta \in \RR$.
\ An infinitely divisible process \ $(X_t)_{t\in\TT}$ \ is said to be
 $(\alpha, \delta)$-dilatively stable if the scaling relation
 \begin{align}\label{dsrel}
  \psi_{Tt_1,\ldots,Tt_k}(\theta_1, \ldots, \theta_k)
  = T^\delta
    \psi_{t_1,\ldots,t_k}(T^{\alpha-\frac{\delta}{2}} \theta_1, \ldots,
                       T^{\alpha-\frac{\delta}{2}} \theta_k)
 \end{align}
 holds for all \ $T > 0$, \ $k \in \NN$,
 \ $\theta_1, \ldots, \theta_k \in \RR$, \ $t_1, \ldots, t_k \in \TT$.
\end{Def}

The functional equation \eqref{dsrel} defines dilative stability in an algebraic sense.
A probabilistic interpretation is given in Remark \ref{Rem_Igl_DS} below.
Note that a process \ $(X_t)_{t\in\TT}$ \ fulfilling the scaling relation
 \eqref{dsrel} with \ $\alpha > 0$ \ and \ $\delta = 0$ \ is an
 \ $\alpha$-self-similar process.
Hence, for infinitely divisible processes, \ $(\alpha, 0)$-dilative stability
 is just \ $\alpha$-self-similarity.
Hence dilative stability is a generalization of self-similarity.
We call the attention that under the conditions of Definition \ref{DefDS}, the
 parameters \ $\alpha$ \ and \ $\delta$ \ are not unique in general, see
 Example \ref{Exl_DSpar_not_unique}, or the simple case of symmetric stable
 processes given below.

To give a very first elementary example for a dilatively stable process, we point out that every L\'evy
 process \ $(L_t)_{t\in\RR_+}$ \ is $(\frac{1}{2}, 1)$-dilatively stable.
Indeed, for all \ $k \in \NN$, \ $0 \leq t_1 < t_2 < \ldots < t_k$ \ and
 \ $T > 0$, \ the increments \ $L_{Tt_j} - L_{Tt_{j-1}}$,
 \ $j \in \{1, \ldots, k\}$, \ are independent, and for each
 \ $j \in \{1, \ldots, k\}$, \ the distribution of \ $L_{Tt_j} - L_{Tt_{j-1}}$
 \ is the $T$-th convolution power of the distribution
 of \ $L_{t_j}-L_{t_{j-1}}$, \ and hence the distribution of
 \ $(L_{Tt_1}, L_{Tt_2} - L_{Tt_1}, \ldots, L_{Tt_k} - L_{Tt_{k-1}})$ \ is the $T$-th
 convolution power of the distribution of
 \ $(L_{t_1}, L_{t_2} - L_{t_1}, \ldots, L_{t_k} - L_{t_{k-1}})$.
\ Consequently, the distribution of \ $(L_{Tt_1}, L_{Tt_2}, \ldots, L_{Tt_k})$
 \ is the $T$-th convolution power of the distribution of
 \ $(L_{t_1}, L_{t_2}, \ldots, L_{t_k})$.
\ In the language of Definition \ref{DefDS}, this means that
 \ $\psi_{Tt_1,\ldots,Tt_k}(\theta_1, \ldots, \theta_k)
    = T \psi_{t_1,\ldots,t_k}(\theta_1, \ldots, \theta_k)$.

In particular, if \ $(L_t)_{t\in\RR_+}$ \ is a symmetric \ $\frac{1}{H}$-stable
 L\'evy process with Hurst index \ $H \geq\frac{1}{2}$, \ then the above
 considerations yield that the process \ $(L_t)_{t\in\RR_+}$ \ is
 \ $(\frac{1}{2}, 1)$-dilatively stable.
Since the process \ $(L_t)_{t\in\RR_+}$ \ is \ $H$-self-similar, it is
 \ $(H, 0)$-dilatively stable as well.
In fact, \ $(L_t)_{t\in\RR_+}$ \ is $(\alpha, \delta)$-dilatively stable for all
 \ $\alpha, \delta \in \RR$ \ satisfying
 \ $\delta + \bigl(\alpha - \frac{\delta}{2}\bigr)\frac{1}{H} = 1$.
\ Indeed, for each \ $t \geq 0$, \ the characteristic exponent of \ $L_t$
 \ has the form \ $\psi_t(\theta) = c t |\theta|^{\frac{1}{H}}$,
 \ $\theta \in \RR$, \ for some \ $c > 0$, \ hence for all \ $T > 0$, \ we
 have
 \ $\psi_{Tt}(\theta) = c T t |\theta|^{\frac{1}{H}}
    = c T^\delta t |T^{\alpha-\frac{\delta}{2}}\theta|^{\frac{1}{H}}
    = T^\delta \psi_t(T^{\alpha-\frac{\delta}{2}} \theta)$,
 \ and independence of the increments yields the scaling relation
 \eqref{dsrel}.
This is an example of a process that is \ $(\alpha, \delta)$-dilatively stable
 with infinitely many choices of \ $(\alpha, \delta)$, \ and \ $\alpha$ \ and
 \ $\delta$ \ can be negative as well.

\begin{Rem}\label{Rem_Igl_DS}
For historical fidelity, we note that Igl\'oi \cite{Igl} originally formulated
 the scaling property \eqref{dsrel} in the language of equality of
 finite-dimensional distributions, namely,
 \begin{align}\label{dsrel_IE}
 \forall\, T>0: \,
  X(T\,\cdot\,) \overset{\text{fd}}
  \sim T^{\alpha-\frac{\delta}{2}} X^{\circledast T^\delta}(\,\cdot\,),
 \end{align}
 where for all \ $c > 0$, \ we denote by \ $X^{\circledast c}$ \ the $c$-th
 convolution power of \ $X$, \ that is, \ $X^{\circledast c}$ \ is a process whose
 finite-dimensional distributions are the \ $c$-th convolution powers
 of the corresponding ones of \ $X$.
\ Comparing the scaling relations \eqref{SS} and \eqref{dsrel_IE}, one can realize that in case of a dilatively stable
 process \ $X$, \ the time scaled process \ $(X_{Tt})_{t\in\TT}$ \ coincides, in the sense of finite dimensional
 distributions, with \ $X$ \ appropriately scaled in space and in convolution exponent as well, whereas in case of
 self-similarity, only some space scaling is allowed.
Roughly speaking, scaling in the convolution exponent is the main new ingredient of dilative stability compared to
 self-similarity.

Igl\'oi \cite{Igl} defined dilative stability under some additional
 assumptions as given in Theorem \ref{Igloi} below.
From this point of view, our Definition \ref{DefDS} is more general.
Note that the condition \ $X_0 = 0$ \ in Theorem \ref{Igloi} turns out to be
 not so restrictive, see Igl\'oi and Barczy \cite[Appendix A]{IglBar} for
 details.
\proofend
\end{Rem}

As it was mentioned, Kaj \cite{Kaj} introduced the notion of aggregate similarity.
Namely, by Definition 1 in Kaj \cite{Kaj}, a stochastic process \ $(X_t)_{t\in\TT}$ \ is
 called aggregate similar with parameter \ $\varrho \in \RR$, \ if it fulfills the scaling relation
 \begin{align}\label{Kaj_scaling}
  \sum_{i=1}^m X^{(i)}(\cdot)
  \overset{\text{fd}}\sim
  m^\varrho X(m^{-\varrho} \, \cdot) \quad\text{ for all }m\in\NN,
 \end{align}
 where \ $(X^{(i)}_t)_{t\in\TT}$, \ $i \in \NN$, \ are independent identically distributed
 copies of \ $(X_t)_{t\in\TT}$.
\ The parameter \ $\varrho$ \ is called rigidity index.
Kaj additionally assumes that \ $\TT=\RR_+$, \ $(X_t)_{t\in\TT}$ \ has stationary
 increments and fulfills \ $X_0 = 0$, \ $\EE(X_t) = 0$ \ for all \ $t \in \TT$,
 \ but these assumptions are not important for our approach.
We introduce the notion of aggregate similarity with parameters \ $\varrho_1,\varrho_2\in\RR$.

\begin{Def}\label{Def_Kaj_aggregate_sim_ext}
Let \ $\varrho_1,\varrho_2\in\RR$.
\ A process \ $(X_t)_{t\in\TT}$ \ is called \ $(\varrho_1,\varrho_2)$-aggregate similar
 if it fulfills the scaling relation
 \begin{align}\label{Kaj_scaling_ext}
  \sum_{i=1}^m X^{(i)}(\cdot)
  \overset{\textup{fd}}\sim
  m^{\varrho_1} X(m^{-\varrho_2} \, \cdot) \quad\text{ for all }m\in\NN,
 \end{align}
 where \ $(X^{(i)}_t)_{t\in\TT}$, \ $i \in \NN$, \ are independent identically distributed
 copies of \ $(X_t)_{t\in\TT}$.
\end{Def}

Under mild regularity assumptions, \eqref{Kaj_scaling_ext}
 already implies the scaling relation \eqref{dsrel} of dilative stability as follows.

\begin{Pro}\label{connection}
\begin{enumerate}
 \item[\textup{(i)}]
  If \ $(X_t)_{t\in\TT}$ \ is a
   \ $(\alpha,\delta)$-dilatively stable
   process for some \ $(\alpha,\delta) \in \RR\times(\RR \setminus \{0\})$, \ then it is
   \ $(\frac12-\frac{\alpha}{\delta},-\frac1{\delta})$-aggregate similar.
 \item[\textup{(ii)}]
  If a process \ $(X_t)_{t\in\TT}$ \ is infinitely divisible, its finite
   dimensional distributions are weakly right-continuous, and it is \ $(\varrho_1,\varrho_2)$-aggregate similar
   for some \ $(\varrho_1,\varrho_2) \in \RR\times(\RR \setminus \{0\})$, \ then it
   is \ $\big(\frac{\varrho_1}{\varrho_2} - \frac{1}{2\varrho_2},-\frac{1}{\varrho_2}\big)$-dilatively
   stable.
\end{enumerate}
\end{Pro}

The proof of Proposition \ref{connection} is given in Section \ref{Proof_lithuania}.
So for infinitely divisible processes, dilative stability and aggregate similarity, given in Definitions
 \ref{DefDS} and \ref{Def_Kaj_aggregate_sim_ext}, are essentially the same.
But note that the purely self-similar case is ruled out in Proposition \ref{connection}, since \ $(\alpha,0)$-dilative stability cannot be handled.

Next we introduce a notion of $(f, g)$-dilative stability in an analogy with
 $f$-self-similarity.

\begin{Def}\label{DefDS_fg}
Let \ $f,g:(0,\infty)\to(0,\infty)$ \ be continuous functions.
An infinitely divisible process \ $(X_t)_{t\in\TT}$ \ is said to be
 $(f,g)$-dilatively stable if the scaling relation
 \begin{align}\label{dsrel_fg}
  \psi_{Tt_1,\ldots,Tt_k}(\theta_1, \ldots, \theta_k)
  = g(T) \psi_{t_1,\ldots,t_k}(f(T) \theta_1, \ldots, f(T) \theta_k)
 \end{align}
 holds for all \ $T > 0$, \ $k \in \NN$,
 \ $\theta_1, \ldots, \theta_k \in \RR$, \ $t_1, \ldots, t_k \in \TT$.
\end{Def}

Theorem 2.2.1 and Proposition 2.1.4 in Igl\'oi \cite{Igl} imply that, under
 some additional assumptions, the functions \ $f$ \ and \ $g$ \ are uniquely
 determined power functions.

\begin{Thm}[Igl\'oi \cite{Igl}]\label{Igloi}
Let \ $(X_t)_{t\in\RR_+}$ \ be an \ $(f, g)$-dilatively stable process.
Suppose that \ $X_0 = 0$, \ $(X_t)_{t\in\RR_+}$ \ is not the identically zero
 process, \ $(X_t)_{t\in\RR_+}$ \ is right continuous in distribution, \ $X_1$
 \ is non-Gaussian and \ $X_t$ \ has finite moments of all orders for all
 \ $t \in \RR_+$, \ depending right continuously on \ $t$.
\ Then there exist unique \ $\alpha \in (0, \infty)$ \ and
 \ $\delta \in [2\alpha, \infty)$ \ such that
 \ $f(T) = T^{\alpha-\frac{\delta}{2}}$,
 \ $T > 0$, \ and \ $g(T) = T^\delta$, \ $T > 0$.
\end{Thm}

Next we give further sufficient conditions under which the functions \ $f$
 \ and \ $g$ \ must be uniquely determined power functions.

\begin{Thm}\label{Thm_DS_fg}
Let \ $(X_t)_{t\in\TT}$ \ be an \ $(f,g)$-dilatively stable process.
\begin{enumerate}
 \item[\textup{(i)}]
  If \ $X$ \ is not the identically zero process and \ $f(T) = T^\beta$ \ for
   some \ $\beta \in \RR$ \ and all \ $T > 0$, \ then there exists a unique
   \ $\gamma \in \RR$ \ such that \ $g(T) = T^\gamma$, \ $T > 0$.
 \item[\textup{(ii)}]
  If \ $X_1$ \ is not Gaussian, but has a Gaussian component in its
   L\'evy--Khintchine representation, then there exist unique
   \ $\beta, \gamma \in \RR$ \ such that \ $f(T) = T^\beta$, \ $T > 0$, \ and
   \ $g(T) = T^\gamma$, \ $T > 0$.
 \item[\textup{(iii)}]
  If \ $\EE(X_1^2) < \infty$, \ $\EE(X_1) \ne 0$ \ and \ $\var(X_1) > 0$,
   \ then there exist unique \ $\beta, \gamma \in \RR$ \ such that
   \ $f(T) = T^\beta$, \ $T > 0$, \ and \ $g(T) = T^\gamma$, \ $T > 0$.
 \item[\textup{(iv)}]
  If \ $X_1$ \ is not Gaussian, \ $\EE(X_1^2) < \infty$ \ and
   \ $\var(X_1) > 0$, \ then there exist unique \ $\beta, \gamma \in \RR$
   \ such that \ $f(T) = T^\beta$, \ $T > 0$, \ and \ $g(T) = T^\gamma$,
   \ $T > 0$.
\end{enumerate}
\end{Thm}

The proof of Theorem \ref{Thm_DS_fg} can be found in Section
 \ref{Proof_Thm_DS_fg}.
Under the assumptions of Theorem \ref{Thm_DS_fg}, the process
 \ $(X_t)_{t\in\TT}$ \ is \ $(T^\beta,T^\gamma)$-dilatively stable
 for some unique \ $\beta,\gamma\in\RR$, \ hence \ $X$ \ could be called
 \ $(\beta,\gamma)$-dilatively stable, as well.
However, for providing an easy link for comparison with Igl\'oi's results
 \cite{Igl}, in Definition \ref{DefDS} we used the parametrization
 \ $(\alpha -\frac{\delta}{2}, \delta)$ \ instead of \ $(\beta, \gamma)$.
Clearly, we have \ $\delta=\gamma$ \ and \ $\alpha=\beta+\frac\delta{2}$.

\begin{Rem}
We give an example of a \ $(f, g)$-dilatively stable process such that \ $f$
 \ and \ $g$ \ are not necessarily power functions.
Namely, let \ $(W_t)_{t\in\RR_+}$ \ be a standard Wiener process.
Then the characteristic exponent of \ $(W_{t_1},\ldots,W_{t_k})$ \ takes the form
 \[
   \psi_{t_1,\ldots,t_k}(\theta_1,\ldots,\theta_k)
   = \frac{1}{2} \sum_{j,\ell=1}^k \min(t_j,t_\ell) \theta_j \theta_\ell
 \]
 for all \ $k \in \NN$, \ $t_1, \ldots, t_k \in \RR_+$, \ and
 \ $\theta_1, \ldots, \theta_k \in \RR$.
\ Hence \ $(W_t)_{t\in\RR_+}$ \ is \ $(f, g)$-dilatively stable for any
 continuous functions \ $f, g : (0, \infty) \to (0, \infty)$ \ satisfying
 \ $g(T) f(T)^2 = T$, \ $T > 0$.
\ Indeed, for all \ $T > 0$,
 \begin{align*}
  \psi_{Tt_1,\ldots,Tt_k}(\theta_1,\ldots,\theta_k)
  &= \frac{1}{2} \sum_{j,\ell=1}^k \min(T t_j, T t_\ell) \theta_j \theta_\ell
   = T \psi_{t_1,\ldots,t_k}(\theta_1, \ldots, \theta_k) \\
  &= g(T) f(T)^2 \psi_{t_1,\ldots,t_k}(\theta_1, \ldots, \theta_k)
   = g(T) \frac{1}{2}
     \sum_{j,\ell=1}^k \min(t_j, t_\ell) f(T) \theta_j\cdot f(T) \theta_\ell\\
  &= g(T) \psi_{t_1,\ldots,t_k}(f(T) \theta_1, \ldots, f(T) \theta_k).
  \hspace*{6.5cm}\proofend
 \end{align*}
\end{Rem}

Next, we recall some dilatively stable processes with stationary increments
 introduced by Igl\'oi \cite{Igl}: fractional L\'evy process, and limits of
 integrated superposition of Ornstein--Uhlenbeck type processes, or of
 continuous-state continuous-time branching processes with immigration (such as
 the Cox--Ingersoll--Ross process).
These processes have the same covariance function as a fractional Brownian
 motion for some parameter $H \in \bigl(\frac{1}{2}, 1\bigr)$.
\ In fact, the non-Gaussian fractional L\'evy process is $(H,1)$-dilatively
 stable, while the other processes are $(H,2H{\,-\,}2)$-dilatively stable.

The remaining part of the paper is organized as follows.
In Section \ref{GFLP} we give new examples for dilatively stable processes
 that are not self-similar, namely, we derive a sufficient condition under
 which a generalized fractional L\'evy process is dilatively stable.
As special cases, sub-fractional L\'evy processes and log-fractional L\'evy
 motions are dilatively stable.
In Section \ref{lithuania} we consider limit processes appearing in certain aggregation
 models under an intermediate growth condition.
In particular, the limit process of joint temporal and contemporaneous aggregation of independent copies of a stationary AR(1) process
 with random-coefficient, introduced in Pilipaysjaut\.{e} and Surgailis \cite{PilSur}, is shown to be dilatively stable.
Section \ref{Proof_lithuania} is devoted to the proof of Proposition \ref{connection}, i.e.,
 the connection between aggregate similarity and dilative stability.
Finally, in Section \ref{Proof_Thm_DS_fg} we prove Theorem \ref{Thm_DS_fg}.

\section{Dilatively stable generalized fractional L\'evy
            processes}
\label{GFLP}

First, we give a summary on two-sided L\'evy processes based on Kl\"uppelberg
 and Matsui \cite[Section 2]{KluMat}.
By a two-sided L\'evy process \ $(L_t)_{t\in\RR}$, \ we mean
 \[
   L_t := L^{(1)}_t \bone_{[0,\infty)}(t) - L^{(2)}_{(-t)-} \bone_{(-\infty,0)}(t) ,
   \qquad t \in \RR ,
 \]
 where \ $L^{(1)}$ \ and \ $L^{(2)}$ \ are independent copies of a L\'evy
 process.
Moreover, we assume that \ $L$ \ is centered (i.e., \ $\EE(L_1) = 0$) without a
 Gaussian component, and its L\'evy measure \ $\mu$ \ satisfies
 \ $\int_{|x|>1} x^2 \, \mu(\dd x) < \infty$, \ hence
  \[
    \EE(L_t^2) = \vert t\vert \EE(L_1^2) =  \vert t\vert \int_\RR x^2 \, \mu(\dd x) < \infty , \qquad
    t \in \RR .
  \]
The above assumptions imply that the characteristic function of \ $L_t$,
 \ $t \in \RR$, \ takes the form
 \ $\EE(\ee^{\ii\theta L_t}) = \ee^{-\vert t\vert\varphi(\theta)}$, \ $t\in\RR$, \ where
 \begin{align}\label{psi}
  \varphi(\theta)
  := \int_\RR (1 - \ee^{\ii\theta x} + \ii \theta x) \, \mu(\dd x) , \qquad
  \theta\in\RR,
 \end{align}
 is the characteristic exponent of \ $L_1$.

The following result is due to Marquardt
 \cite[Proposition 2.1 and (2.13)]{Mar}.

\begin{Pro}\label{Pro_GFLP}
Let \ $L=(L_t)_{t\in\RR}$ \ be a two-sided L\'evy process without a Gaussian
 component.
Assume that \ $\EE(L_1) = 0$ \ and \ $\EE(L_1^2) < \infty$.
\ Let \ $f : \RR^2 \to \RR$ \ be a measurable function such that
 \ $\RR \ni u \mapsto f(t,u) \in L^2(\RR)$ \ for all \ $t \in \RR$.
\ Then the integral \ $S_t := \int_\RR f(t,u) \, L(\dd u)$ \ exists in
 \ $L^2(\Omega, \cF, \PP)$-sense for all \ $t \in \RR$.
\ Furthermore, for \ $s, t \in \RR$ \ we obtain \ $\EE(S_t) = 0$, \ the
 isometry
 \[
   \EE(S_t^2) = \EE(L_1^2) \int_\RR f(t,u)^2 \,\dd u, \qquad t \in \RR ,
 \]
 and
 \[
   \cov(S_s, S_t) = \EE(L_1^2) \int_\RR f(s,u) f(t,u) \, \dd u, \qquad
   s, t \in \RR .
 \]
Moreover, the characteristic function of \ $(S_{t_1}, \ldots, S_{t_k})$,
 \ $t_1 < t_2 < \cdots < t_k$, \ $t_j \in \RR$, \ $j \in \{1, \ldots, k\}$,
 \ takes the form
 \begin{align}\label{GFLP_char}
  \EE\left(\exp\left\{\ii \sum_{j=1}^k \theta_j S_{t_j}\right\}\right)
  = \exp\left\{-\int_\RR
                 \varphi\left(\sum_{j=1}^k \theta_j f(t_j,u)\right)
                 \dd u\right\}
 \end{align}
 for \ $\theta_j \in \RR$, \ $j \in \{1, \ldots, k\}$, \ where
 \ $\varphi$ \ is given by \eqref{psi}.
\end{Pro}

\begin{Def}
A stochastic process \ $(S_t)_{t\in\RR}$ \ given in Proposition \ref{Pro_GFLP}
 is called a generalized fractional L\'evy process corresponding to the kernel
 function \ $f$.
\end{Def}

We note that Kl\"uppelberg and Matsui \cite[Definition 2.2]{KluMat} have
 already introduced generalized fractional L\'evy processes with kernel
 function \ $f(t,u) := g(t-u) - g(-u)$, \ $t, u \in \RR$, \ for some function
 \ $g : \RR \to \RR$ \ satisfying \ $g(t) = 0$ \ for \ $t < 0$ \ and
 \ $\int_\RR (g(t-u) - g(-u))^2 \, \dd u < \infty$.

Next we derive a sufficient condition under which a generalized fractional
 L\'evy process is dilatively stable.

\begin{Pro}\label{condDS}
If the kernel function \ $f$ \ given in Proposition \ref{Pro_GFLP} satisfies
 \begin{align}\label{DS_funct_eq}
  f(t,u)
  = T^{\alpha-\frac{\delta}{2}} f\left(\frac{t}{T}, \frac{u}{T^\delta}\right) ,
  \qquad t,u\in\RR,\quad T>0,
 \end{align}
 then the generalized fractional L\'evy process \ $(S_t)_{t\in\RR}$ \ is
 \ $(\alpha, \delta)$-dilatively stable.
If, in addition, \ $\EE(L_1^2) > 0$ \ and
 \ $\int_\RR f(1, u)^2 \, \dd u > 0$, \ then the parameters of dilative
 stability  \ $\alpha$ \ and \ $\delta$ \ are uniquely determined.
\end{Pro}

\noindent{\bf Proof.}
First observe that the process \ $S$ \ is infinitely divisible, since by
 \eqref{GFLP_char}, for all \ $k \in \NN$,
 \ $\theta_1, \ldots, \theta_k \in \RR$, \ $t_1, \ldots, t_k \in \RR$ \ and
 \ $n \in \NN$, \ we have
 \[
   \EE\left(\exp\left\{\ii \sum_{j=1}^k \theta_j S_{t_j}\right\}\right)
   = \left(\exp\left\{-\frac{1}{n}
                       \int_\RR
                        \varphi\left(\sum_{j=1}^k \theta_j f(t_j,u)\right)
                        \dd u\right\}\right)^n ,
 \]
 where, again by \eqref{GFLP_char},
 \begin{align*}
  \exp\left\{ -\frac{1}{n}
               \int_\RR
                \varphi\left(\sum_{j=1}^k \theta_j f(t_j,u)\right)
                \dd u\right\}
  &= \exp\left\{-\int_\RR
                  \varphi\left(\sum_{j=1}^k \theta_j f(t_j, nv)\right)
                  \dd v\right\} \\
  &= \EE\left(\exp\left\{\ii \sum_{j=1}^k \theta_j \tS_{t_j}\right\}\right)
 \end{align*}
 is the characteristic function of \ $(\tS_{t_1}, \ldots, \tS_{t_k})$ \ with
 \ $\tS_t:= \int_\RR f(t,nv) \, L(\dd v)$, \ $t \in \RR$.

Again by \eqref{GFLP_char}, the characteristic exponent of
 \ $(S_{t_1}, \ldots, S_{t_k})$ \ has the form
 \[
   \psi_{t_1,\ldots,t_k}(\theta_1, \ldots, \theta_k)
   = \int_\RR \varphi\left(\sum_{j=1}^k \theta_j f(t_j, u)\right) \dd u ,
   \qquad \theta_1, \ldots, \theta_k \in \RR .
 \]
Hence for all \ $T > 0$, \ $k \in \NN$, \ $\theta_1, \ldots, \theta_k \in \RR$
 \ and \ $t_1, \ldots, t_k \in \RR$, \ by \eqref{DS_funct_eq}, we get
 \begin{align*}
  &T^\delta
   \psi_{t_1,\ldots,t_k}(T^{\alpha-\frac{\delta}{2}} \theta_1, \ldots,
                      T^{\alpha-\frac{\delta}{2}} \theta_k)
   = T^\delta
      \int_\RR
       \varphi\left(\sum_{j=1}^k T^{\alpha-\frac{\delta}{2}} \theta_j f(t_j, u)\right)
       \dd u \\
  &\qquad
   = \int_\RR
      \varphi\left(\sum_{j=1}^k
                    T^{\alpha-\frac{\delta}{2}} \theta_j
                    f(Tt_j/T, v/T^\delta)\right)
      \dd v
   = \int_\RR \varphi\left(\sum_{j=1}^k \theta_j f(Tt_j, v)\right) \dd v \\
  &\qquad
   = \psi_{Tt_1,\ldots,Tt_k}(\theta_1, \ldots, \theta_k) ,
 \end{align*}
 which yields that \ $(S_t)_{t\in\RR}$ \ is \ $(\alpha, \delta)$-dilatively
 stable.
If, in addition, \ $\EE(L_1^2)>0$ \ and \ $\int_\RR f(1,u)^2\,\dd u>0$, \ then,
 by Proposition \ref{Pro_GFLP}, \ $\EE(S_1^2)\in(0,\infty)$.
\ Using also that \ $S_1$ \ is not Gaussian, part (iv) of Theorem
 \ref{Thm_DS_fg} yields the uniqueness of \ $\alpha$ \ and \ $\delta$.
\proofend

Next we formulate three examples for generalized fractional L\'evy processes
 that are dilatively stable, but not self-similar.

\begin{Ex}({\bf Sub-fractional L\'evy process})
Sub-fractional Brownian motions were introduced by Bojdecki et al.\
 \cite{BojGorTal}, see also Tudor \cite[Section 1.3]{Tud}.
Let \ $(B^{(H)}_t)_{t\in\RR}$ \ be a fractional Brownian motion on \ $\RR$ \ with
 parameter \ $H \in \bigl(\frac{1}{2}, 1\bigr)$, \ i.e.,
 \ $(B^{(H)}_t)_{t\in\RR}$ \ is a Gaussian process with zero mean having
 covariance function
 \[
   \cov(B^{(H)}_t , B^{(H)}_s)
   = \frac{1}{2} \bigl(|t|^{2H} + |s|^{2H} - |t - s|^{2H} \bigr) ,
   \qquad t, s \in \RR .
 \]
Let
 \[
   S^{(H)}_t := \frac{1}{\sqrt{2}} (B^{(H)}_t + B^{(H)}_{-t}) , \qquad
   t \in \RR_+ .
 \]
Then \ $(S^{(H)}_t)_{t\in\RR_+}$ \ is a sub-fractional Brownian motion with
 parameter \ $H$, \ i.e., it is a Gaussian process with zero mean and with
 covariance function
 \[
   \cov(S^{(H)}_t , S^{(H)}_s)
   = t^{2H} + s^{2H} - \frac{1}{2} \bigl((s+t)^{2H} + |t - s|^{2H}\bigr) \qquad
   t, s \in \RR_+ ,
 \]
 see, e.g., Tudor \cite[Proposition 1.17]{Tud}.
The process \ $S^{(H)}$ \ is $H$-self-similar, see, e.g., Tudor
 \cite[Proposition 1.14]{Tud}.
Recall that
 \[
   B^{(H)}_t
   = \frac{\sqrt{2H+1} \sin(\pi H)}{\Gamma\bigl(H+\frac{1}{2}\bigr)}
     \int_\RR \bigl((t-u)_+^{H-\frac{1}{2}} - (-u)_+^{H-\frac{1}{2}}\bigr) B(\dd u) ,
   \qquad t \in \RR ,
 \]
 where \ $B$ \ is a Brownian motion on \ $\RR$, \ see Samorodnitsky and Taqqu
 \cite[Proposition 7.2.6]{SamTaq} and Beran et al.\
 \cite[Section 1.3.5]{BerFenGhoKul}.
Hence we get
 \begin{align*}
  S^{(H)}_t
  = \frac{\sqrt{2H+1} \sin(\pi H)}
         {\sqrt{2}\, \Gamma\bigl(H+\frac{1}{2}\bigr)}
    \biggl[ &\int_\RR
              \bigl((t-u)_+^{H-\frac{1}{2}} - (-u)_+^{H-\frac{1}{2}}\bigr) B(\dd u)\\
          &+\int_\RR
             \bigl((-t-u)_+^{H-\frac{1}{2}} - (-u)_+^{H-\frac{1}{2}}\bigr) B(\dd u)
    \biggr], \qquad t \in \RR_+ .
 \end{align*}
Motivated by this moving average representation of \ $S^{(H)}$, \ we introduce
 so-called sub-fractional L\'evy processes.
Namely, for all \ $H \in \bigl(\frac{1}{2}, 1\bigr)$, \ let
 \begin{align*}
  \widetilde S^{(H)}_t
  := \frac{1}{\sqrt{2}} \bigl(L_t^{(H)} + L_{-t}^{(H)}\bigr), \qquad t \in \RR_+ ,
 \end{align*}
 where
 \[
   L_t^{(H)} := \frac{\sqrt{2H+1} \sin(\pi H)}{\Gamma\bigl(H+\frac{1}{2}\bigr)}
               \int_\RR
                \bigl((t-u)_+^{H-\frac{1}{2}} - (-u)_+^{H-\frac{1}{2}}\bigr)
                L(\dd u) ,
   \qquad t \in \RR ,
 \]
 with a two-sided L\'evy process  \ $(L_t)_{t\in\RR}$ \ without a Gaussian
 component satisfying \ $\EE(L_1) = 0$ \ and \ $0 < \EE(L_1^2) < \infty$.
\ Note that \ $L_t^{(H)}$ \ is well-defined, since the kernel function
 \ $f : \RR^2 \to \RR$,
 \[
   f(t,u):=
   \begin{cases}
    (t-u)_+^{H-\frac{1}{2}} - (-u)_+^{H-\frac{1}{2}} ,
    & \text{if \ $t \ne u$ \ and \ $u \ne 0$,} \\
    0 , & \text{if \ $t = u$ \ or \ $u = 0$,}
   \end{cases}
 \]
 satisfies the conditions of Proposition \ref{Pro_GFLP}, see Beran et al.\
 \cite[Section 3.7.2.5]{BerFenGhoKul}.

We call the process \ $\widetilde S^{(H)}$ \ a sub-fractional L\'evy process
 with parameter \ $H$ \ corresponding to the two-sided L\'evy process \ $L$.
\ The process \ $\widetilde S^{(H)}$ \ is \ $(H, 1)$-dilatively stable by
 Proposition \ref{condDS}, and the parameters are uniquely determined, since
 \ $\EE(L_1^2) > 0$, \ and the kernel function \ $f$ \ satisfies condition
 \eqref{DS_funct_eq} with \ $\alpha = H$ \ and \ $\delta = 1$ \ together with
 \ $\int_\RR f(1,u)^2 \, \dd u > 0$.
\ We note that the process \ $\widetilde S^{(H)}$ \ has long memory increments,
 since \ $H \in \bigl(\frac{1}{2}, 1\bigr)$, \ see Kl\"uppelberg and Matsui
 \cite[Lemma 2.6]{KluMat}.
\proofend
\end{Ex}

\begin{Ex}({\bf Log-fractional L\'evy motion})
Let \ $\alpha\in(1,2)$ \ and \ $M$ \ be a symmetric \ $\alpha$-stable random measure with Lebesgue control measure.
The process
 \[
   X_t:=\int_\RR \bigl(\log(|t - u|) - \log(|u|)\bigr) M(\dd u) , \qquad
   t \in \RR ,
 \]
 is called a symmetric log-fractional stable motion, see Samorodnitsky and
 Taqqu \cite[Example 3.6.6]{SamTaq}.
For the definition of the integral with respect to \ $M$ \ see
 \cite[Section 3.4]{SamTaq}, and the fact that \ $X_t$, \ $t\in\RR$, \ is
 well defined can be checked in the same way as in Samorodnitsky and Taqqu
 \cite[Example 3.6.5]{SamTaq}.
Especially,
 \ $\int_\RR \bigl(\log(|t - u|) - \log(|u|)\bigr)^2 \, \dd u < \infty$,
 \ $t \in \RR$.
\ The process \ $X$ \ is \ $\frac{1}{\alpha}$-self-similar with stationary
 increments.
Motivated by this, let us introduce so-called log-fractional L\'evy motions.
Namely, we call the process
 \[
   \tX_t := \int_\RR \bigl(\log(|t - u|) - \log(|u|)\bigr) L(\dd u) , \qquad
   t \in \RR ,
 \]
 a log-fractional L\'evy motion, where \ $(L_t)_{t\in\RR}$ \ is a two-sided
 L\'evy process without a Gaussian component satisfying \ $\EE(L_1) = 0$ \ and
 \ $0 < \EE(L_1^2) <\infty$.
\ Note that \ $\tX_t$ \ is well-defined, since the kernel function
 \ $f : \RR^2 \to \RR$,
 \[
   f(t,u):=
   \begin{cases}
    \log(|t - u|) - \log(|u|) , & \text{if \ $t \ne u$ \ and \ $u \ne 0$,} \\
    0 , & \text{if \ $t = u$ \ or \ $u = 0$,}
   \end{cases}
 \]
 satisfies the conditions of Proposition \ref{Pro_GFLP}, see Samorodnitsky and
 Taqqu \cite[Example 3.6.6]{SamTaq}.

The process \ $(\tX_t)_{t\in\RR_+}$ \ is
 \ $\bigl(\frac{1}{2}, 1\bigr)$-dilatively stable by Proposition \ref{condDS},
 and the parameters are uniquely determined, since \ $\EE(L_1^2) > 0$, \ and,
 for the kernel function \ $f$, \ condition \eqref{DS_funct_eq} holds with
 \ $\alpha = \frac{1}{2}$ \ and \ $\delta = 1$ \ together with
 \ $\int_\RR f(1,u)^2 \, \dd u > 0$.
\proofend
\end{Ex}

\begin{Ex}
Generalized sub-fractional Brownian motions were introduced by Sghir
 \cite{Sgh}.
In one of the representation theorems of Sghir \cite[Theorem 2.2]{Sgh}, a
 self-similar process comes into play that can be generalized to present
 another example for a dilatively stable generalized fractional L\'evy process
 as follows.
For \ $K \in (0, 2)$, \ let
 \begin{align*}%\label{X}
  X^{(K)}_t := \int_{\RR_+} (1 - \ee^{-ut}) u^{-\frac{K+1}{2}} \, B(\dd u) , \qquad
  t \in \RR_+ ,
 \end{align*}
 where \ $(B_t)_{t\in\RR_+}$ \ is a standard Brownian motion.
The process \ $ X^{(K)}$ \ is $\frac{K}{2}$-self-similar.
Let
 \begin{align*}%\label{tX}
  \tX^{(K)}_t := \int_{\RR_+} (1 - \ee^{-ut}) u^{-\frac{K+1}{2}} \, L(\dd u) , \qquad
  t \in \RR_+ ,
 \end{align*}
 where \ $(L_t)_{t\in\RR_+}$ \ is a two-sided L\'evy process without a Gaussian
 component satisfying \ $\EE(L_1) = 0$ \ and \ $0 < \EE(L_1^2)<\infty$.
\ Note that \ $\tX^{(K)}_t$ \ is well-defined, since the kernel function
 \ $f : \RR^2 \to \RR$,
 \[
   f(t,u)
   :=\begin{cases}
      (1 - \ee^{-ut}) u^{-\frac{K+1}{2}} , & \text{if \ $t > 0$ \ and \ $u > 0$,}\\
      0 , & \text{if \ $t \leq 0$ \ or \ $u \leq 0$,}
     \end{cases}
 \]
 satisfies the conditions of Proposition \ref{Pro_GFLP}.
Indeed, we must prove that
 \[
   \int_{\RR_+} (1 - \ee^{-ut})^2 u^{-(K+1)} \, \dd u < \infty , \qquad t > 0 .
 \]
This integral may diverge at \ $u = 0$ \ or at \ $u = \infty$.
\ As \ $u \to \infty$, \ the integrand behaves like \ $u^{-(K+1)}$, \ which is
 integrable, since \ $\int_1^\infty u^{-(K+1)} \, \dd u < \infty$.
\ As \ $u \downarrow 0$, \ the integrand behaves like
 \ $u^2 u^{-(K+1)} = u^{-K+1}$, \ which is integrable, since
 \ $\int_0^1 u^{-K+1} \, \dd u < \infty$.

The process \ $\tX^{(K)}$ \ is \ $\bigl(\frac{K}{2}, -1\bigr)$-dilatively
 stable by Proposition \ref{condDS}, and the parameters are uniquely
 determined, since \ $\EE(L_1^2) > 0$, \ and one can check that condition
 \eqref{DS_funct_eq} holds for the kernel function \ $f$ \ with
 \ $\alpha = \frac{K}{2}$ \ and \ $\delta = -1$ \ together with
 \ $\int_\RR f(1,u)^2 \, \dd u > 0$.
\proofend
\end{Ex}

In the next example we point out that under the conditions of Definition
 \ref{DefDS}, the parameters \ $\alpha$ \ and \ $\delta$ \ are not unique in
 general.

\begin{Ex}\label{Exl_DSpar_not_unique}
Let \ $\alpha \in (0, 2]$ \ and \ $M$ \ be a symmetric \ $\alpha$-stable random measure with Lebesgue control measure.
For \ $H \in (0, 1)$,  \ $H \ne \frac{1}{\alpha}$, \ a well-balanced linear
 fractional stable process \ $X=(X_t)_{t\in\RR}$\ is given by
 \[
   X_t := \int_\RR
           \big(|t-u|^{H-\frac{1}{\alpha}} - |u|^{H-\frac{1}{\alpha}}\big) M(\dd u) ,
   \qquad t \in \RR ,
 \]
 see, e.g., Samorodnitsky and Taqqu \cite[Example 3.6.5]{SamTaq}.
On the one hand, it is known that \ $X$ \ is \ $H$-self-similar and thus
 \ $(H,0)$-dilatively stable, since the process \ $X$ \ is infinitely
 divisible.
Further, the characteristic exponent of \ $(X_{t_1}, \ldots, X_{t_k})$ \ can be
 written in the form
 \[
   \psi_{t_1,\ldots,t_k}(\theta_1, \ldots, \theta_k)
   = \int_\RR \left|\sigma \sum_{j=1}^k \theta_j f(t_j, u)\right|^\alpha \dd u
 \]
 for \ $t_1 < \cdots < t_k$ \ and \ $\theta_1, \ldots, \theta_k \in \RR$ \ for
 some \ $\sigma > 0$, \ see, e.g., Samorodnitsky and Taqqu
 \cite[Property 3.2.1]{SamTaq}, where the kernel function
 \ $f : \RR^2 \to \RR$, \ is given by
 \[
   f(t,u):=
   \begin{cases}
    |t-u|^{H-1/\alpha} - |u|^{H-1/\alpha} ,
    & \text{if \ $t \ne u$ \ and \ $u \ne 0$,} \\
    0 , & \text{if \ $t = u$ \ or \ $u = 0$.}
   \end{cases}
 \]
Hence for all \ $T > 0$, \ $k \in \NN$,
 \ $\theta_1, \ldots, \theta_k \in \RR$, \ and \ $t_1,\ldots,t_k\in\RR$, \ we
 get
 \begin{align*}
  \psi_{Tt_1,\ldots,Tt_k}(\theta_1, \ldots, \theta_k)
  & = \int_\RR
       \left|\sigma \sum_{j=1}^k \theta_j f(T t_j, u)\right|^\alpha \dd u
    = T \int_\RR
       \left|\sigma \sum_{j=1}^k \theta_j f(T t_j, T v)\right|^\alpha \dd v \\
  & = T \int_\RR
         \left|\sigma \sum_{j=1}^k T^{H-1/\alpha} \theta_j f(t_j, v) \right|^\alpha
         \dd v
    = T
      \psi_{t_1,\ldots,t_k}(T^{H-1/\alpha} \theta_1, \ldots, T^{H-1/\alpha} \theta_k) ,
 \end{align*}
 showing that \ $X$ \ is
 \ $\bigl(H - \frac{1}{\alpha} + \frac{1}{2}, 1\bigr)$-dilatively stable.
This example shows that the parametrization of dilatively stable processes is
 not unique in general.
However, note also that the process \ $M$ \ does not have finite second
 moments for \ $\alpha<2$ \  and that for $\alpha=2$ the linear fractional
 Brownian motion \ $X$ \ is Gaussian.
\proofend
\end{Ex}

\section{Further examples from aggregation models}\label{lithuania}

In this section we consider certain processes arising by aggregation of independent copies of random systems, appropriately rescaled in space and time. The behavior of such systems can vary significantly depending on the relative speed of aggregation with respect to the rescaling. We are particularly interested in intermediate regimes, where the speed of aggregation and rescaling is in a certain balance.
In this way non self-similar, but dilatively stable limit processes may appear as follows, serving as further examples.

Our first example comes from a joint temporal and contemporaneous aggregation
\begin{equation}\label{jtca}
S_{N,n}(t)=\sum_{i=1}^N\sum_{s=1}^{\lfloor nt\rfloor}X_i(s) , \qquad t\geq0,
\end{equation}
of \ $N$ \ independent and identically distributed, stationary random coefficient AR(1) processes
 \ $(X_i(s))_{s\in\RR_+}$, $i\in\{1,\ldots,N\}$, \ with some given mixing density,
 recently investigated by Pilipaysjaut\.{e} and Surgailis \cite{PilSur}. Under an intermediate growth condition
\begin{equation}\label{interm}
\frac{N^{1/(1+\beta)}}{n}\to c\in(0,\infty)
\end{equation}
as simultaneously \ $n,N\to\infty$, \ by Theorem 2.2 in \cite{PilSur} for any \ $\beta\in(-1,1)$ \ the appropriately rescaled system \eqref{jtca} converges to a random process \ $(\sqrt{c}Z_\beta(t/c))_{t\in\RR_+}$ \ in the sense of finite-dimensional distributions.
Here the process \ $Z_\beta$ \ is infinitely divisible by Proposition 3.1 in \cite{PilSur} and determined by the characteristic
 exponent of \ $(Z_\beta(t_1), \ldots, Z_\beta(t_k))$ \ which takes the form
 \begin{align*}
   \psi_{t_1,\ldots,t_k}(\theta_1,\ldots,\theta_k)
      = - C \int_0^\infty \!\!\left(\! \exp\left\{\!-\frac{1}{2}\int_\RR \left(\sum_{j=1}^k \theta_j (f(x,t_j-s) - f(x,-s))
           \right)^2\dd s \right\} - 1\right)x^\beta\,\dd x,
 \end{align*}
 for all \ $t_j\in\RR_+$, \ $\theta_j\in\RR$, \ $j=1,\ldots,k$, \ $k\in\NN$, and some constant \ $C>0$, \
  where \ $f:\RR^2\to\RR$ \ is given by
 \begin{equation}\label{fdef}
     f(x,t):=\begin{cases}
             \frac{1-\ee^{-xt}}{x} & \text{if \ $x>0$ \ and \ $t>0$,}\\
             0                     & \text{otherwise.}
            \end{cases}
 \end{equation}
 As remarked on page 1022 of \cite{PilSur} the process \ $Z_\beta$ \ is not self-similar and not stable.

\begin{Pro}\label{Pro_lithuanian}
For any \ $\beta\in(-1,1)$, \ the process \ $Z_\beta$ \ is \ $\left(1-\frac{\beta}{2},-\beta-1\right)$-dilatively stable.
\end{Pro}

\noindent{\bf Proof.}
By Proposition 3.1 in Pilipaysjaut\.{e} and Surgailis \cite{PilSur}, \ $Z_\beta$ \ is infinitely divisible.
Note that the function \ $f$ \ in \eqref{fdef} fulfills \ $f(x,Tt)=T\cdot f(Tx,t)$ \ for all $x\in\RR$, $t\in\RR$ and $T>0$. Hence,
for all \ $T>0$, \ $k\in\NN$, \ $\theta_1,\ldots,\theta_k\in\RR$ \ and \ $t_1,\ldots,t_k\in\RR_+$, \
 by change of variables \ $s=Tu$ \ and then \ $Tx=y$, \ we get
 \begin{align*}
  &\psi_{Tt_1,\ldots,Tt_k}(\theta_1,\ldots,\theta_k)\\
     &  = - C \int_0^\infty \Bigg( \exp\Bigg\{-\frac{1}{2}\int_{\RR} \left(\sum_{j=1}^k \theta_j
                                       \big(f(x,Tt_j-s)-f(x,-s)\big)\right)^2 \dd s
                                        \Bigg\} - 1\Bigg)x^\beta\,\dd x\\
     &  = - C \int_0^\infty \Bigg( \exp\Bigg\{-\frac{1}{2}\int_{\RR} \left(\sum_{j=1}^k T^{3/2}\theta_j
                                       \big(f(Tx,t_j-u)-f(Tx,-u)\big)\right)^2 \dd u
                                        \Bigg\} - 1\Bigg)x^\beta\,\dd x\\
     & = - T^{-\beta-1}C \int_0^\infty \Bigg( \exp\Bigg\{-\frac{1}{2}\int_{\RR} \left(\sum_{j=1}^k T^{3/2}\theta_j
                                       \big(f(y,t_j-u)-f(y,-u)\big)\right)^2 \dd u
                                        \Bigg\} - 1\Bigg)y^\beta\,\dd y\\
     & = T^{-\beta-1} \psi_{t_1,\ldots,t_k}(T^{3/2}\theta_1,\ldots,T^{3/2}\theta_k).
 \end{align*}
 By definition of dilative stability, this yields the assertion.
\proofend

Note that, by Proposition 3.1 in \cite{PilSur}, for \ $\beta\in(0,1)$ \ the process \ $Z_\beta$ \ is centered with
 \ $\EE(Z_\beta(1)^2)\in(0,\infty)$, \ and for  \ $\beta\in(-\frac12,1)$ \  a.s.\ has continuous trajectories.
Since \ $Z_\beta(1)$ \ is not Gaussian, part (iv) of Theorem \ref{Thm_DS_fg} shows that at least for \ $\beta\in(0,1)$ \ the parameters of dilative stability for \ $Z_\beta$ \ are unique.

An alternative way to verify dilative stability of \ $Z_\beta$ \ is by means of aggregate similarity introduced
in Definition \ref{Def_Kaj_aggregate_sim_ext}.
As remarked on page 1023 in \cite{PilSur}, the process \ $(Z_\beta(t^{2/3}))_{t\geq0}$ \ is
 \ $(\frac{3}{2(1+\beta)},\frac{3}{2(1+\beta)})$-aggregate similar, where \ $\beta\in(-1,1)$.
\ Hence, by an easy calculation, \ $Z_\beta$ \ is \ $(\frac{3}{2(1+\beta)},\frac{1}{1+\beta})$-aggregate similar.
 Then, since \ $Z_\beta$ \ is infinitely divisible (see Proposition 3.1 in \cite{PilSur}),
 part (ii) of Proposition \ref{connection} shows that \ $Z_\beta$ \ is
 \ $\left(1-\frac{\beta}{2},-\beta-1\right)$-dilatively stable at least for the reduced parameter set
 \ $\beta\in(-\frac12,1)$, \ where weak right-continuity is justified by a.s.\ continuity of the sample paths.

Note that, by Theorems 2.1 and 2.2 in \cite{PilSur}, the limit process of the appropriately rescaled system \eqref{jtca}
 under slow and fast growth condition, i.e., \ $c=0$ \ or \ $c=\infty$ \ in \eqref{interm}, or
 under iterated aggregation \ $\lim_{n\to\infty}\lim_{N\to\infty}$ \ or \ $\lim_{N\to\infty}\lim_{n\to\infty}$ \ is either fractional Brownian motion, or a linear time multiple of a stable random variable, or a variance mixture of Brownian motion with a stable mixing variable. Since all of these processes are self-similar, they are naturally dilatively stable as well.
All in all, each of the limit processes in Theorems 2.1 and 2.2 in \cite{PilSur} is dilatively stable.

Kaj's motivation to introduce the notion of aggregate similarity in \cite{Kaj} comes from an aggregation model having its origin
 in the study of total workload in teletraffic models.
The survey in \cite{Kaj} shows that in several different models (superposition of renewal counting processes,
 sums of inverse L\'evy subordinators, infinite source Poisson models, self-similar rate models, inference model for wireless communication) under an intermediate growth condition the same scaling limit process \ $Y_\gamma$ \ arises, depending on a stability parameter \ $\gamma\in(1,2)$; \ cf.\ also \cite{Gai,GaiKaj,KajML,KajTaq}.
By Theorem 2 in \cite{GaiKaj},
 this process is particularly not self-similar and not stable but \ $Y_\gamma$ \ has a.s.\ continuous trajectories.
 Since \ $Y_\gamma$ \ is aggregate similar with rigidity index \ $\varrho=(\gamma-1)^{-1}$ \ by Section 3.6 in \cite{Kaj}
 and is infinitely divisible by Section 3.4 in \cite{Gai}, it is also $(\frac{3-\gamma}{2},1-\gamma)$-dilatively stable
  by part (ii) of Proposition \ref{connection}.
Further, Gaigalas \cite[Sections 3.1, 3.2 and 4]{Gai}
 has shown that the process \ $Y_\gamma$ \ builds a certain bridge between a $\gamma$-stable L\'evy process and fractional Brownian motion with Hurst index \ $H=\frac{3-\gamma}{2}$ \ appearing as scaling limits in the above mentioned models under slow, respectively fast growth condition.
Note that the $\gamma$-stable L\'evy process is also $(\frac{3-\gamma}{2},1-\gamma)$-dilatively stable (see the paragraph before Remark
 \ref{Rem_Igl_DS}) and the $\frac{3-\gamma}{2}$-selfsimilar fractional Brownian motion is
 also $(\frac{3-\gamma}{2},0)$-dilatively stable.
Hence all the possible limit processes are $(\alpha,\delta)$-dilatively stable with the same parameter
 $\alpha = \frac{3-\gamma}{2}$ giving a justification for Igl\'oi's parametrization (see the comment after Theorem \ref{Thm_DS_fg}).

Since, by Theorem 4 in \cite{Kaj}, the process \ $Y_\gamma$ \ has an integral representation with respect to a compensated Poisson random measure for which a natural extension to a multivariate spatial model has been given  by Bierm\'e et al. \cite{BieEstKaj}, it seems to be possible to generalize the notion of dilative stability to higher dimensions. The characteristic function of the random field in question can be found on page 1139 in \cite{BieEstKaj} and suggests that the scaling relation of dilative stability is fulfilled. However, we renounce to consider this question in full detail, since for random fields even a more general scaling by linear operators as in \cite{BMS} might be introduced, giving rise for future research.

\section{Proof of Proposition \ref{connection}}
\label{Proof_lithuania}

First, note that in terms of characteristic exponents, for an infinitely
 divisible process \ $(X_t)_{t\in\TT}$, \ the scaling relation
 \eqref{Kaj_scaling_ext} is equivalent to
 \begin{align}\label{Kaj_scaling_char_exp}
  m \, \psi_{t_1,\ldots,t_k}(\theta_1, \ldots, \theta_k)
  = \psi_{m^{-\varrho_2}t_1,\ldots,m^{-\varrho_2}t_k}(m^{\varrho_1} \theta_1, \ldots,
                                          m^{\varrho_1} \theta_k)
 \end{align}
 for every \ $m \in \NN$, \ $k \in \NN$,
 \ $\theta_1, \ldots, \theta_k \in \RR$, \ and \ $t_1, \ldots, t_k \in \TT$.

Suppose that \ $(X_t)_{t\in\TT}$ \ is a
 \ $(\alpha,\delta)$-dilatively stable
 process for some \  $(\alpha,\delta) \in \RR \times (\RR \setminus \{0\})$.
\ By \eqref{dsrel},
 \begin{align*}
  \psi_{Tt_1,\ldots,Tt_k}(\theta_1, \ldots, \theta_k)
  = T^{\delta}
    \psi_{t_1,\ldots,t_k}(T^{\alpha-\frac{\delta}{2}} \theta_1, \ldots, T^{\alpha-\frac{\delta}{2}} \theta_k)
 \end{align*}
 for all \ $T > 0$, \ $k \in \NN$, \ $\theta_1, \ldots, \theta_k \in \RR$,
 \ and \ $t_1, \ldots, t_k \in \TT$.
\ By choosing \ $T = m^{1/\delta}$, \ $m \in \NN$, \ we have
 \[
   \psi_{m^{1/\delta}t_1,\ldots,m^{1/\delta}t_k}(\theta_1, \ldots, \theta_k)
   = m \,\psi_{t_1,\ldots,t_k}(m^{\frac{\alpha}{\delta}-\frac12} \theta_1, \ldots, m^{\frac{\alpha}{\delta}-\frac12} \theta_k) .
 \]
Replacing \ $\theta_i$ \ by \ $m^{\frac12-\frac{\alpha}{\delta}} \theta_i$,
 \ $i \in \{1, \ldots, k\}$, \ we have \eqref{Kaj_scaling_char_exp} with \ $\varrho_1=\frac12-\frac{\alpha}{\delta}$ \ and \ $\varrho_2=-\frac1{\delta}$, \ concluding part (i).

Let us now suppose that \ $(X_t)_{t\in\TT}$ \ is an infinitely divisible, weakly
 right-continuous
 \ $(\varrho_1, \varrho_2)$-aggregate similar process with some
 \ $(\varrho_1,\varrho_2) \in \RR \times (\RR \setminus \{0\})$.
\ Then, by letting \ $s_i := m^{-\varrho_2} t_i$,
 \ $\widetilde\theta_i := m^{\varrho_1} \theta_i$, \ $i \in \{1, \ldots, k\}$, \ we
 have that \eqref{Kaj_scaling_char_exp} is equivalent to
 \begin{align}\label{Kaj_scaling_char_exp2}
  \psi_{m^{\varrho_2} s_1,\ldots,m^{\varrho_2} s_k}(m^{-\varrho_1} \widetilde\theta_1, \ldots,
                                    m^{-\varrho_1} \widetilde\theta_k)
  = m^{-1} \psi_{s_1,\ldots,s_k}(\widetilde\theta_1, \ldots, \widetilde\theta_k)
 \end{align}
 for all \ $m \in \NN$, \ $k \in \NN$, \ $s_1, \ldots, s_k \in \TT$ \ and
 \ $\widetilde\theta_1, \ldots, \widetilde\theta_k \in \RR$.
\ Now let \ $z = \frac{m}{n}$ \ be an arbitrary positive rational number.
By \eqref{Kaj_scaling_char_exp} and \eqref{Kaj_scaling_char_exp2}, we have
 \begin{align}\label{Kaj_scaling_char_exp3}
  \begin{split}
   z \,\psi_{t_1,\ldots,t_k}(\theta_1, \ldots, \theta_k)
   &= n^{-1} \psi_{m^{-\varrho_2}t_1,\ldots,m^{-\varrho_2}t_k}(m^{\varrho_1} \theta_1, \ldots,
                                                  m^{\varrho_1} \theta_k) \\
   &= n^{-1}
      \psi_{n^{-\varrho_2}z^{-\varrho_2}t_1,\ldots,n^{-\varrho_2}z^{-\varrho_2}t_k}
       (n^{\varrho_1} z^{\varrho_1} \theta_1, \ldots, n^{\varrho_1} z^{\varrho_1} \theta_k) \\
   &= \psi_{z^{-\varrho_2}t_1,\ldots,z^{-\varrho_2}t_k}(z^{\varrho_1} \theta_1, \ldots,
                                            z^{\varrho_1}\theta_k)
  \end{split}
 \end{align}
 for all \ $k \in \NN$, \ $t_1, \ldots, t_k \in \TT$ \ and
 \ $\theta_1, \ldots, \theta_k \in \RR$.

Finally, if \ $z > 0$ \ is arbitrary, then choose a sequence \ $(z_n)_{n\in\NN}$
 \ of positive rational numbers such that \ $z_n \downarrow z$ \ as
 \ $n \to \infty$.
\ By \eqref{Kaj_scaling_char_exp3}, we get
 \begin{align*}
  z \, \psi_{t_1,\ldots,t_k}(\theta_1, \ldots, \theta_k)
  &= \lim_{n\to\infty} z_n \psi_{t_1,\ldots,t_k}(\theta_1, \ldots, \theta_k)
   = \lim_{n\to\infty}
      \psi_{z_n^{-\varrho_2}t_1,\ldots,z_n^{-\varrho_2}t_k}(z_n^{\varrho_1} \theta_1, \ldots,
                                               z_n^{\varrho_1}\theta_k) \\
  &= \psi_{z^{-\varrho_2}t_1,\ldots,z^{-\varrho_2}t_k}(z^{\varrho_1} \theta_1, \ldots,
                                          z^{\varrho_1} \theta_k) ,
 \end{align*}
 where the last equality can be checked as follows.
By L\'evy's continuity theorem, weak right-continuity of the finite-dimensional
 distributions of \ $(X_t)_{t\in\TT}$ \ yields
 \ $\EE(\exp\{\ii \sum_{j=1}^k \theta_j X_{z_n^{-\varrho_2}t_j}\})
    \to \EE(\exp\{\ii \sum_{j=1}^k \theta_j X_{z^{-\varrho_2}t_j}\})$
 \ as \ $n \to \infty$ \  for all \ $k \in \NN$,
 \ $t_1, \ldots, t_k \in \TT$, \ uniformly on compact subsets of
 \ $(\theta_1, \ldots, \theta_k)^\top \in \RR^k$.
\ Further, by Lemma 3.2.11 in Stroock \cite{Str}, we conclude continuity of
 the characteristic exponent \ $\psi_{t_1,\ldots,t_k}$ \ for all \ $k \in \NN$
 \ and \ $t_1, \ldots, t_k \in \TT$, \ and
 \ $\psi_{z_n^{-\varrho_2}t_1,\ldots,z_n^{-\varrho_2}t_k}(\theta_1, \ldots, \theta_k)
    \to \psi_{z^{-\varrho_2}t_1,\ldots,z^{-\varrho_2}t_k}(\theta_1, \ldots, \theta_k)$
 \ as \ $n \to \infty$ \ uniformly on compacts for all \ $k \in \NN$,
 \ $t_1, \ldots, t_k \in \TT$ \ and \ $\theta_1, \ldots, \theta_k \in \RR$.
Consequently, denoting the closed ball around
 \ $(z^{\varrho_1} \theta_1, \ldots, z^{\varrho_1}\theta_k)$ \ with radius \ $1$ \ by
 \ $B^{z,\varrho_1}_{\theta_1,\ldots,\theta_k}$, \ we have, for sufficiently large
 \ $n \in \NN$,
 \begin{align*}
  &|\psi_{z_n^{-\varrho_2}t_1,\ldots,z_n^{-\varrho_2}t_k}
     (z_n^{\varrho_1} \theta_1, \ldots, z_n^{\varrho_1}\theta_k)
    - \psi_{z^{-\varrho_2}t_1,\ldots,z^{-\varrho_2}t_k}
       (z^{\varrho_1} \theta_1, \ldots, z^{\varrho_1} \theta_k)| \\
  &\qquad\qquad
   \leq |\psi_{z_n^{-\varrho_2}t_1,\ldots,z_n^{-\varrho_2}t_k}
          (z_n^{\varrho_1} \theta_1, \ldots, z_n^{\varrho_1}\theta_k)
         - \psi_{z^{-\varrho_2}t_1,\ldots,z^{-\varrho_2}t_k}
            (z_n^{\varrho_1} \theta_1, \ldots, z_n^{\varrho_1} \theta_k)| \\
  &\qquad\qquad
   \quad
        + |\psi_{z^{-\varrho_2}t_1,\ldots,z^{-\varrho_2}t_k}
            (z_n^{\varrho_1} \theta_1, \ldots, z_n^{\varrho_1}\theta_k)
           - \psi_{z^{-\varrho_2}t_1,\ldots,z^{-\varrho_2}t_k}
              (z^{\varrho_1} \theta_1, \ldots, z^{\varrho_1} \theta_k)|
   \end{align*}
  \begin{align*}
  &\qquad\qquad
   \leq \sup_{x\in B^{z,\varrho_1}_{\theta_1,\ldots,\theta_k}}
         |\psi_{z_n^{-\varrho_2}t_1,\ldots,z_n^{-\varrho_2}t_k}(x)
          - \psi_{z^{-\varrho_2}t_1,\ldots,z^{-\varrho_2}t_k}(x)| \\
 &\qquad\qquad
   \quad
        + |\psi_{z^{-\varrho_2}t_1,\ldots,z^{-\varrho_2}t_k}
            (z_n^{\varrho_1} \theta_1, \ldots, z_n^{\varrho_1}\theta_k)
           - \psi_{z^{-\varrho_2}t_1,\ldots,z^{-\varrho_2}t_k}
              (z^{\varrho_1} \theta_1, \ldots, z^{\varrho_1} \theta_k)| \\
 &\qquad\qquad
  \to 0 \qquad \text{as \ $n \to \infty$.}
 \end{align*}
Letting \ $T = z^{-\varrho_2}$ \ and \ $\widetilde\theta_i = z^{\varrho_1} \theta_i$,
 \ $i \in \{1, \ldots, k\}$, \ we have
 \[
   \psi_{Tt_1,\ldots,Tt_k}(\widetilde\theta_1, \ldots, \widetilde\theta_k)
   = T^{-\frac{1}{\varrho_2}}
     \psi_{t_1,\ldots,t_k}(T^{\frac{\varrho_1}{\varrho_2}} \widetilde\theta_1, \ldots, T^{\frac{\varrho_1}{\varrho_2}} \widetilde\theta_k)
 \]
 for all \ $k \in \NN$, \ $T > 0$, \ $t_1, \ldots, t_k \in \TT$ \ and
 \ $\widetilde\theta_1, \ldots, \widetilde\theta_k \in \RR$, \ which coincides
 with the definition of
 \ $\bigl(\frac{\varrho_1}{\varrho_2} - \frac{1}{2\varrho_2}, -\frac{1}{\varrho_2}\bigr)$-dilative stability
 concluding part (ii).
\proofend

\section{Proof of Theorem \ref{Thm_DS_fg}}\label{Proof_Thm_DS_fg}

Applying \eqref{dsrel_fg} with \ $k = 1$ \ and \ $t_1 = 1$, \ we get
 \begin{align}\label{help1}
  \psi_T(\theta) = g(T) \psi_1(f(T) \theta) ,
  \qquad \theta \in \RR , \quad T \in (0, \infty) .
 \end{align}
Replacing \ $T$ \ by \ $ts$ \ with \ $s,t\in(0,\infty)$, \ we have
 \begin{align*}
  \psi_{ts}(\theta) = g(ts) \psi_1(f(ts) \theta) ,
  \qquad \theta \in \RR , \quad s, t \in (0, \infty) .
 \end{align*}
Using \eqref{dsrel_fg} with \ $k = 1$, $T = t$ \ and $t_1 = s$, \ and applying
 \eqref{help1}, we get
 \begin{align*}
  \psi_{ts}(\theta)
  = g(t) \psi_s(f(t) \theta)
  = g(t) g(s) \psi_1(f(t) f(s) \theta) ,
  \qquad \theta \in \RR , \quad s, t \in (0, \infty) .
 \end{align*}
The uniqueness of the characteristic exponent of the infinitely divisible
 random variable \ $X_{ts}$ \ implies that
 \begin{equation}\label{gpsif}
  g(st) \psi_1(f(st) \theta) = g(s)g(t) \psi_1(f(s) f(t) \theta),
  \qquad \theta \in \RR , \quad s, t \in (0, \infty) .
 \end{equation}

If there exists a \ $\beta \in \RR$ \ such that \ $f(T) = T^\beta$,
 \ $T \in(0, \infty)$, \ then \eqref{gpsif} takes the form
 \begin{equation}\label{gpsi}
  g(st) \psi_1((s t)^\beta \theta) = g(s) g(t) \psi_1((s t)^\beta \theta) ,
  \qquad \theta \in \RR , \quad s, t \in (0, \infty) .
 \end{equation}
If \ $X$ \ is not the identically zero process, then \ $\psi_1$ \ is not the
 identically zero function, hence there exists \ $\theta_0 \in \RR$ \ such
 that \ $\psi_1(\theta_0) \ne 0$.
\ Substituting \ $\theta = \theta_0/(st)^\beta$ \ into \eqref{gpsi} and dividing
 both sides by \ $\psi_1(\theta_0)$, \ we obtain \ $g(st) = g(s) g(t)$ \ for
 all \ $s, t \in (0, \infty)$.
\ Since \ $g$ \ is supposed to be continuous, it is known that there exists
 some \ $\gamma \in \RR$ \ such that \ $g(t) = t^\gamma$, \ $t \in (0, \infty)$.
\ Substituting \ $\theta = \theta_0/T^\beta$ \ into \eqref{help1} and dividing
 both sides by \ $\psi_1(\theta_0)$, \ we get
 \ $g(T) = \psi_T(\theta_0/T^\beta)/\psi_1(\theta_0)$ \ for all
 \ $T \in (0, \infty)$, \ hence we obtain the uniqueness of \ $\gamma$, \ and
 we conclude statement (i).

Recall that there exist some \ $a \in \RR$, \ $\sigma \in \RR_{++}$ \ and a
 measure \ $\mu$ \ on \ $\RR \setminus \{0\}$ \ satisfying
 \ $\int_{\RR\setminus\{0\}} \min(1, x^2) \, \mu(\dd x) < \infty$ \ (called L\'evy
 measure) such that
 \begin{align}\label{help_char_exp}
  \psi_1(\theta)
  = \ii a\theta + \frac{1}{2}\sigma^2 \theta^2
    + \int_{\RR\setminus\{0\}}
       \left(1 - \ee^{\ii\theta x} + \frac{\ii\theta x}{1+x^2}\right) \mu(\dd x) ,
  \qquad \theta \in \RR .
 \end{align}
Hence we get
 \begin{align*}
  g(st) \psi_1(f(s t) \theta)
  &=\ii g(st) f(st)a \theta + \frac{1}{2} g(st) f(st)^2 \sigma^2 \theta^2 \\
  &\quad
    + g(st)
      \int_{\RR\setminus\{0\}}
       \left(1 - \ee^{\ii f(st)\theta x} + \frac{\ii f(st)\theta x}{1+x^2}\right)
       \mu(\dd x) \\
  &=\ii \widetilde a \theta
    + \frac{1}{2} \widetilde\sigma^2 \theta^2
    + \int_{\RR\setminus\{0\}}
       \left(1 - \ee^{\ii\theta x} + \frac{\ii\theta x}{1+x^2}\right)
       \widetilde\mu(\dd x) ,
    \qquad \theta \in \RR , \quad s, t \in (0, \infty) ,
 \end{align*}
 where \ $\widetilde\sigma^2 := g(st) f(st)^2 \sigma^2$,
 \ $\widetilde\mu(\dd x) := g(st) \mu\left(\frac{\dd x}{f(st)}\right)$ \ is a
 L\'evy measure, and
 \begin{align*}
  \ii \widetilde a \theta
  &:= \ii g(st)f(st)a\theta
      + g(st)
        \int_{\RR\setminus\{0\}}
         \left(1 - \ee^{\ii f(st)\theta x} + \frac{\ii f(st)\theta x}{1+x^2}\right)
         \mu(\dd x) \\
  &\quad
      - \int_{\RR\setminus\{0\}}
         \left(1 - \ee^{\ii\theta x} + \frac{\ii\theta x}{1+x^2}\right)
         \widetilde\mu(\dd x) \\
  &=\ii \theta g(s t)
    \left(f(s t) a
          + \int_{\RR\setminus\{0\}}
             \left(\frac{y}{1+(y/f(st))^2} - \frac{y}{1+y^2} \right)
             \mu\left(\frac{\dd y}{f(st)}\right)\right) \\
  &=\ii \theta g(st) f(st)
    \left(a + \int_{\RR\setminus\{0\}}
               \left(\frac{x}{1+x^2} - \frac{x}{1+f(st)^2x^2} \right)
               \mu(\dd x)\right) .
 \end{align*}
Similarly,
 \begin{align*}
  g(s) g(t) \psi_1(f(s) f(t) \theta)
  = \ii \widehat a \theta
    + \frac{1}{2} \widehat\sigma^2 \theta^2
    + \int_{\RR\setminus\{0\}}
       \left(1 - \ee^{\ii\theta x} + \frac{\ii \theta x}{1+x^2}\right)
       \widehat\mu(\dd x) ,
 \end{align*}
 for all \ $\theta \in \RR$ \ and \ $s, t \in (0, \infty)$, \ where
 \ $\widehat\sigma^2 := g(s) g(t) f(s)^2 f(t)^2 \sigma^2$,
 \ $\widehat\mu(\dd x) = g(s) g(t) \mu\left(\frac{\dd x}{f(s)f(t)}\right)$
 \ is a L\'evy measure, and
 \begin{align*}
  \widehat a
  := g(s) g(t) f(s) f(t)
     \left(a + \int_{\RR\setminus\{0\}}
                \left(\frac{x}{1+x^2} - \frac{x}{1+f(s)^2f(t)^2x^2} \right)
                \mu(\dd x)\right) .
 \end{align*}
Since the triplet corresponding to the infinitely divisible random variable
 \ $X_{st}$ \ is uniquely determined provided that the truncation function is
 fixed (see, e.g., Stroock \cite[Theorem 3.2.20]{Str}), we have
 \ $\widetilde a = \widehat a$, \ $\widetilde\sigma^2 = \widehat\sigma^2$ \ and
 \ $\widetilde\mu = \widehat\mu$.

Now we suppose that \ $X_1$ \ is not Gaussian, but has a Gaussian component in
 its L\'evy--Khintchine representation, i.e., \ $\sigma \ne 0$ \ and
 \ $\mu \ne 0$.
\ Then \ $\widetilde\sigma^2 = \widehat\sigma^2$ \ yields
 that \ $h(st) = h(s) h(t)$, \ $s, t \in (0, \infty)$, \ for the function
 \ $h : (0,\infty) \to (0,\infty)$, \ $h(t) := g(t) f(t)^2$,
 \ $t \in (0, \infty)$.
\ Since \ $f$ \ and \ $g$ \ are supposed to be continuous, it is known that
 there exists some \ $\kappa \in \RR$ \ such that \ $g(t) f(t)^2 = t^\kappa$,
 \ $t \in (0, \infty)$.
\ Then \ $\widetilde\mu = \widehat\mu$ \ implies
 \ $g(st) \mu\left(\frac{\dd x}{f(st)}\right)
    = g(s) g(t) \mu\left(\frac{\dd x}{f(s)f(t)}\right)$ \ for all
 \ $s, t \in (0, \infty)$.
\ Since \ $\mu$ \ is a L\'evy measure, we have
 \ $\mu\bigl((y, \infty)\bigr) < \infty$ \ for all \ $y \in (0, \infty)$,
 \ hence
 \[
   g(s t) \mu\left(\left(\frac{y}{f(s t)}, \infty\right)\right)
   = g(s) g(t) \mu\left(\left(\frac{y}{f(s)f(t)}, \infty\right)\right) ,
   \qquad y, s, t \in (0, \infty) .
 \]
Consequently, \ $\mu\bigl((z, \infty)\bigr) = a \mu\bigl((bz, \infty)\bigr)$
 \ for all \ $z \in (0, \infty)$, \ where
 \ $a := g(s)g(t)/g(st) \in (0, \infty)$ \ and
 \ $b := f(st)/(f(s)f(t)) \in (0, \infty)$.
\ We have \ $a = b^2$, \ since we have already proven that
 \ $g(t) f(t)^2 = t^\kappa$ \ for all \ $t \in (0, \infty)$.
\ The aim of the following discussion is to show that \ $a \in (0, 1)$ \ would
 lead to a contradiction.
If so, then clearly \ $b \in (0, 1)$.
\ Iterating \ $\mu(\dd z) = a \mu(b \dd z)$, \ we obtain
 \ $\mu(\dd z) = a^n \mu(b^n \dd z)$ \ for all \ $n \in \NN$.
\ Since \ $\mu \ne 0$, \ without loss of generality, we may suppose that
 \ $\mu\bigl((0, \infty)\bigr) > 0$.
\ Then there exists \ $z_0 \in (0, \infty)$ \ such that
 \ $\mu\bigl((z_0, \infty)\bigr) > 0$.
\ Then \ $\mu\bigl((z_0, \infty)\bigr) = a^n \mu\bigl((b^n z_0, \infty)\bigr)$
 \ for all \ $n \in \NN$.
\ Since \ $\mu\bigl((z_0, \infty)\bigr) > 0$, \ $a^n \to 0$ \ and
 \ $\mu\bigl((b^n z_0, \infty)\bigr) \to \mu\bigl((0, \infty)\bigr)$ \ as
 \ $n \to \infty$, \ we conclude \ $\mu\bigl((0, \infty)\bigr) = \infty$.
\ Consequently, using that \ $\mu$ \ is a L\'evy measure,
 \ $\int_{(0,1]} z^2 \, \mu(\dd z) \in (0, \infty)$.
\ Applying \ $\mu(\dd z) = a^n \mu(b^n \dd z)$, \ $n \in \NN$, \ we obtain
 \begin{align*}
  \int_{(0,1]} z^2 \, \mu(\dd z)
  &= \int_{(0,1]} z^2 \, a^n \mu(b^n \dd z)
  = a^n \int_{(0,b^n]} \left(\frac{y}{b^n}\right)^2 \, \mu(\dd y) \\
  &= \left(\frac{a}{b^2}\right)^n \int_{(0,b^n]} y^2 \, \mu(\dd y)
   = \int_{(0,b^n]} y^2 \, \mu(\dd y) ,
 \end{align*}
 and \ $\lim_{n\to\infty} \int_{(0,b^n]} y^2 \, \mu(\dd y) = 0$ \ leads to a
 contradiction.
In a similar way, \ $a \in (1, \infty)$ \ would also lead to a contradiction,
 since then we would have
 \ $\mu\bigl((z, \infty)\bigr) = a^{-1} \mu\bigl((b^{-1}z, \infty)\bigr)$ \ for
 all \ $z \in (0, \infty)$, \ where \ $a^{-1} \in (0, 1)$.
\ Summarizing, the only possibility is \ $a = 1$, \ which implies \ $b = 1$,
 \ and hence \ $f(st) = f(s) f(t)$, \ $s, t \in (0, \infty)$, \ and
 \ $g(st) = g(s) g(t)$, \ $s, t \in (0, \infty)$.
\ Since \ $f$ \ and \ $g$ \ are supposed to be continuous, it is known that
 there exist some \ $\beta, \gamma \in \RR$ \ such that \ $f(t) = t^\beta$,
 \ $t \in (0, \infty)$, \ and \ $g(t) = t^\gamma$, \ $t \in (0, \infty)$.
\ Now, we show that \ $\beta$ \ and \ $\gamma$ \ are uniquely determined.
Using \eqref{help1} and \eqref{help_char_exp}, we obtain
 \begin{align*}
  \psi_T(\theta)
  &= \ii a T^{\beta+\gamma} \theta
     + \frac{1}{2} \sigma^2 T^{2\beta+\gamma} \theta^2
     + T^\gamma
       \int_{\RR\setminus\{0\}}
        \left(1-\ee^{\ii T^\beta\theta x} + \frac{\ii T^\beta\theta x}{1+x^2}\right)
        \mu(\dd x) \\
  &= \ii \left( a
               + \int_{\RR\setminus\{0\}}
                  \left(\frac{x}{1+x^2} - \frac{x}{1+T^{2\beta}x^2}\right)
                  \mu(\dd x) \right)
     T^{\beta+\gamma} \theta
     + \frac{1}{2} \sigma^2 T^{2\beta+\gamma} \theta^2 \\
  &\quad
     + \int_{\RR\setminus\{0\}}
        \left(1-\ee^{\ii\theta x} + \frac{\ii\theta x}{1+x^2}\right)
        T^\gamma \mu\left(\frac{\dd x}{T^\beta}\right)
 \end{align*}
 for all \ $\theta \in \RR$ \ and \ $T \in (0, \infty)$.
\ If there exist some \ $\tbeta, \tgamma \in \RR$ \ such that the process
 \ $(X_t)_{t\in\TT}$ \ is \ $(\tf, \tg)$-dilatively stable with
 \ $\tf(t) = t^{\tbeta}$, \ $t \in (0, \infty)$, \ and
 \ $g(t) = t^{\tgamma}$, \ $t \in (0, \infty)$, \ then, by the same reasoning,
 \begin{align*}
  \psi_T(\theta)
  &= \ii \left( a
               + \int_{\RR\setminus\{0\}}
                  \left(\frac{x}{1+x^2} - \frac{x}{1+T^{2\tbeta}x^2}\right)
                 \mu(\dd x) \right)
     T^{\tbeta+\tgamma} \theta
     + \frac{1}{2} \sigma^2 T^{2\tbeta+\tgamma} \theta^2 \\
  &\quad
     + \int_{\RR\setminus\{0\}}
        \left(1-\ee^{\ii\theta x} + \frac{\ii\theta x}{1+x^2}\right)
        T^{\tgamma} \mu\left(\frac{\dd x}{T^{\tbeta}}\right)
 \end{align*}
 for all \ $\theta \in \RR$ \ and \ $T \in (0, \infty)$.
\ Since the triplet corresponding to the infinitely divisible random variable
 \ $X_T$ \ is uniquely determined provided that the truncation function is
 fixed (see, e.g., Stroock \cite[Theorem 3.2.20]{Str}), we have
 \ $\sigma^2 T^{2\beta+\gamma} = \sigma^2 T^{2\tbeta+\tgamma}$, \ $T \in (0, \infty)$,
 \ and
 \ $T^\gamma \mu\left(\frac{\dd x}{T^\beta}\right)
    = T^{\tgamma} \mu\left(\frac{\dd x}{T^{\tbeta}}\right)$,
 \ $T \in (0, \infty)$.
\ Thus \ $\sigma \ne 0$ \ implies \ $2\beta + \gamma = 2\tbeta + \tgamma$,
 \ and hence \ $\tgamma - \gamma = 2(\beta - \tbeta)$, \ which yields
 \ $\mu(\dd y) = T^{\tgamma-\gamma} \mu(T^{\beta-\tbeta} \dd y)
    = T^{2(\beta-\tbeta)} \mu(T^{\beta-\tbeta} \dd y)$, \ $T \in (0, \infty)$.
\ The aim of the following discussion is to show that \ $\beta < \tbeta$
 \ would lead us to a contradiction.
Since \ $\mu \ne 0$, \ without loss of generality, we may suppose that
 \ $\mu\bigl((0, \infty)\bigr) > 0$.
\ Then there exists \ $y_0 \in (0, \infty)$ \ such that
 \ $\mu\bigl((y_0, \infty)\bigr) > 0$.
\ Then
 \ $\mu\bigl((y_0, \infty)\bigr)
    = T^{2(\beta-\tbeta)} \mu\bigl((T^{\beta-\tbeta} y_0, \infty)\bigr)$
 \ for all \ $T \in (0, \infty)$.
\ Since \ $\mu\bigl((y_0, \infty)\bigr) > 0$, \ $T^{2(\beta-\tbeta)} \to 0$ \ and
 \ $\mu\bigl((T^{\beta-\tbeta} y_0, \infty)\bigr) \to \mu\bigl((0, \infty)\bigr)$
 \ as \ $T \to \infty$, \ we conclude \ $\mu\bigl((0, \infty)\bigr) = \infty$.
\ Consequently, using again that \ $\mu$ \ is a L\'evy measure,
 \ $\int_{(0,1]} y^2 \, \mu(\dd y) \in (0, \infty)$.
\ Applying \ $\mu(\dd y) = T^{2(\beta-\tbeta)} \mu(T^{\beta-\tbeta} \dd y)$,
 \ $T \in (0, \infty)$, \ we obtain
 \begin{align*}
  \int_{(0,1]} y^2 \, \mu(\dd y)
  &= \int_{(0,1]} y^2 \, T^{2(\beta-\tbeta)} \mu(T^{\beta-\tbeta} \dd y)
  = T^{2(\beta-\tbeta)}
    \int_{(0,T^{\beta-\tbeta}]}
     \left(\frac{z}{T^{\beta-\tbeta}}\right)^2 \, \mu(\dd z) \\
  &= \int_{(0,T^{\beta-\tbeta}]} z^2 \, \mu(\dd z) ,
 \end{align*}
 and \ $\lim_{T\to\infty} \int_{(0,T^{\beta-\tbeta}]} y^2 \, \mu(\dd y) = 0$ \ leads
 us to a contradiction.
In a similar way, \ $\beta > \tbeta$ \ would also lead us to a contradiction,
 since then we would have
 \ $\mu\bigl((z, \infty)\bigr)
    = T^{2(\tbeta-\beta)} \mu\bigl((T^{\tbeta-\beta}z, \infty)\bigr)$
 \ for all \ $z \in (0, \infty)$.
\ Summarizing, the only possibility is \ $\beta = \tbeta$, \ which implies
 \ $\gamma = \tgamma$, \ and we conclude statement (ii).

Next we suppose that \ $\EE(X_1^2) < \infty$, \ $\EE(X_1) \ne 0$ \ and
 \ $\var(X_1) > 0$.
\ Then \ $\psi_1$ \ is twice differentiable and \ $\EE(X_T) = \ii \psi_T'(0)$,
 \ $T \in (0, \infty)$, \ and \ $\var(X_T) = \psi_T''(0)$,
 \ $T \in (0, \infty)$.
\ Taking the derivative of \eqref{gpsif} and putting \ $\theta = 0$, \ we
 obtain
 \[
   g(st) f(st) \psi_1'(0) = g(s) g(t) f(s) f(t) \psi_1'(0) ,
   \quad s, t \in (0, \infty) .
 \]
Since \ $\psi_1'(0) = -\ii \EE(X_1) \ne 0$, \ we conclude that there exists
 some \ $\lambda \in \RR$ \ such that \ $g(t) f(t) = t^\lambda$,
 \ $t \in (0, \infty)$.
\ Taking the second derivative of \eqref{gpsif} and putting \ $\theta = 0$,
 \ we get
 \[
   g(st) f(st)^2 \psi_1''(0) = g(s) g(t) f(s)^2 f(t)^2 \psi_1''(0) ,
   \quad s, t \in (0, \infty) .
 \]
Since \ $\psi_1''(0) = \var(X_1) \ne 0$, \ we obtain that there exists some
 \ $\kappa \in \RR$ \ such that \ $g(t) f(t)^2 = t^\kappa$,
 \ $t \in (0, \infty)$.
\ Since \ $f$ \ and \ $g$ \ are positive, it yields that
 \ $f(t) = t^{\kappa - \lambda}$, \ $t \in (0, \infty)$, \ and
 \ $g(t) = t^{2\lambda-\kappa}$, \ $t \in (0, \infty)$.
\ Choosing \ $\beta := \kappa - \lambda$ \ and
 \ $\gamma := 2 \lambda - \kappa$, \ we have \ $f(t) = t^\beta$,
 \ $t \in (0, \infty)$, \ and \ $g(t) = t^\gamma$, \ $t \in (0, \infty)$.
\ Next we show that \ $\beta$ \ and \ $\gamma$ \ are uniquely determined.
Using \eqref{help1}, \ we get
 \begin{align}\label{gpsibetagamma}
  \psi_T(\theta) = T^\gamma \psi_1(T^\beta \theta) ,
  \qquad \theta \in \RR , \qquad T \in (0, \infty) .
 \end{align}
Taking the derivative of \eqref{gpsibetagamma} and putting \ $\theta = 0$,
 \ we obtain \ $\psi_T'(0) = T^{\beta+\gamma} \psi_1'(0)$, \ $T \in (0, \infty)$.
\ Since \ $\psi_1'(0) = -\ii \EE(X_1) \ne 0$, \ we conclude that
 \ $\beta + \gamma$ \ is uniquely determined.
Taking the second derivative of \eqref{gpsibetagamma} and putting
 \ $\theta = 0$, \ we obtain \ $\psi_T''(0) = T^{2\beta+\gamma} \psi_1''(0)$,
 \ $T \in (0, \infty)$.
\ Since \ $\psi_1''(0) = \var(X_1) \ne 0$, \ we have that \ $2\beta + \gamma$
 \ is also uniquely determined, thus we get the uniqueness of \ $\beta$
 \ and \ $\gamma$, \ and we conclude statement (iii).

Finally, we suppose that \ $X_1$ \ is not Gaussian, \ $\EE(X_1^2) < \infty$
 \ and \ $\var(X_1) > 0$.
\ As in the proof of (iii), there exists some \ $\kappa \in \RR$
 \ such that \ $g(t) f(t)^2 = t^\kappa$, \ $t \in (0, \infty)$;
 \ further, provided that there exist some \ $\beta, \gamma \in \RR$ \ such
 that \ $f(t) = t^\beta$, \ $t \in (0, \infty)$, \ and \ $g(t) = t^\gamma$,
 \ $t \in (0, \infty)$, \ then \ $2 \beta + \gamma$ \ is uniquely determined.
The proof of statement (iv) can be finished as in the case of (ii).
\proofend

\begin{Rem}
{(i)}
The proof of part (iii) of Theorem \ref{Thm_DS_fg} uses the same ideas as that
 of Theorem 2.2.1 in Igl\'oi \cite{Igl}, however, instead of cumulants we work
 with moments, and we have much weaker (moment) assumptions.

\noindent {(ii)}
Similarly to the proof of part (ii) of Theorem \ref{Thm_DS_fg}, one can check
 that if \ $(X_t)_{t\in\TT}$ \ is an $(\alpha,\delta)$-dilatively stable process
 for some \ $\alpha, \delta \in\RR$ \ such that \ $\EE(X_1^2) < \infty$,
 \ then
 \[
   \EE(X_{Tt}) = T^{\alpha+\frac{\delta}{2}} \EE(X_t) \qquad \text{and} \qquad
   \var(X_{Tt}) = T^{2\alpha} \var(X_t), \qquad T > 0 , \quad t \in \TT .
 \]
\end{Rem}

\section*{Acknowledgements}
We are grateful to the referee for several valuable comments and suggestions, in particular
 for calling our attention to the examples presented in Section \ref{lithuania}.

\end{document}